\documentclass{gtart}
\usepackage{amssymb,amsmath}

\input gtoutput
\volumenumber{2}\papernumber{6}\volumeyear{1998}
\pagenumbers{103}{116}\published{12 July 1998}
\proposed{Dieter Kotschick}\seconded{Tomasz Mrowka, John Morgan}
\received{27 February 1998}
\accepted{9 July 1998}

\newtheorem{thm}{Theorem}[section]
\newtheorem{cor}[thm]{Corollary}
\newtheorem{lem}[thm]{Lemma}
\newtheorem{prop}[thm]{Proposition}

\newtheorem{rem}[thm]{Remark}

\newtheorem{quest}[thm]{Question}
\newtheorem{conj}[thm]{Conjecture}
\numberwithin{equation}{section}

\newcommand{\ga}{\gamma}
\newcommand{\Ga}{\Gamma}

\newcommand{\de}{\delta}
\newcommand{\ep}{\epsilon}

\newcommand{\La}{\Lambda}

\newcommand{\Om}{\Omega}
\newcommand{\om}{\omega}
\newcommand{\ro}{\rho}
\newcommand{\si}{\sigma}

\renewcommand{\th}{\theta}

\newcommand{\ze}{\zeta}

\renewcommand{\t}{\mathbf t}
\newcommand{\s}{\mathbf s}

\newcommand{\C}{\mathbb C}
\newcommand{\Z}{\mathbb Z}
\newcommand{\Q}{\mathbb Q}
\newcommand{\R}{\mathbb R}
\newcommand{\x}{\times}

\newcommand{\del}{\partial}

\newcommand{\Hom}{\operatorname{\mathrm Hom}} 

\newcommand{\rank}{\operatorname{\mathrm rank}}

\newcommand{\SW}{\operatorname{\mathrm SW}}

\newcommand{\Tor}{\operatorname{\mathrm Tor}}

\newcommand{\Spin}{\operatorname{\mathrm Spin}}

\newcommand{\Id}{\operatorname {\mathrm Id}}
\newcommand{\End}{\operatorname {\mathrm End}}

\begin{document}
%

\title{Symplectic fillings and positive scalar curvature}

\author{Paolo Lisca} 

\address{Dipartimento di Matematica\\ Universit\`a di Pisa \\I-56127
Pisa, ITALY} 

\email{lisca@dm.unipi.it}

\asciiaddress{Dipartimento di Matematica\\ Universita di Pisa \\I-56127
Pisa, ITALY}

\asciiabstract{%
Let X be a 4-manifold with contact boundary. We prove that the
monopole invariants of X introduced by Kronheimer and Mrowka vanish
under the following assumptions: (i) a connected component of the
boundary of X carries a metric with positive scalar curvature and
(ii) either b_2^+(X)>0 or the boundary of X is disconnected. As an
application we show that the Poincare homology 3-sphere, oriented
as the boundary of the positive E_8 plumbing, does not carry
symplectically semi-fillable contact structures. This proves, in
particular, a conjecture of Gompf, and provides the first example of a
3-manifold which is not symplectically semi-fillable. Using work
of Froyshov, we also prove a result constraining the topology of
symplectic fillings of rational homology 3-spheres having positive
scalar curvature metrics.}

\begin{abstract} 
Let $X$ be a $4$--manifold with contact boundary. We prove that the
monopole invariants of $X$ introduced by Kronheimer and Mrowka vanish
under the following assumptions: (i) a connected component of the
boundary of $X$ carries a metric with positive scalar curvature and
(ii) either $b_2^+(X)>0$ or the boundary of $X$ is disconnected. As an
application we show that the Poincar\'e homology $3$--sphere, oriented
as the boundary of the positive $E_8$ plumbing, does not carry
symplectically semi-fillable contact structures. This proves, in
particular, a conjecture of Gompf, and provides the first example of a
$3$--manifold which is not symplectically semi-fillable. Using work
of Fr\o yshov, we also prove a result constraining the topology of
symplectic fillings of rational homology $3$--spheres having positive
scalar curvature metrics.
\end{abstract}

\primaryclass{53C15}
\secondaryclass{57M50, 57R57} 

\keywords{Contact structures, monopole equations, Seiberg--Witten
equations, positive scalar curvature, symplectic fillings}

\asciikeywords{Contact structures, monopole equations, Seiberg-Witten
equations, positive scalar curvature, symplectic fillings}

\maketitlepage
%
%

\section{Introduction}\label{s:intro}

\subsection{Basic facts and questions on contact structures}\label{ss:contact}

Let $Y$ be a closed $3$--manifold. A coorientable field of $2$--planes
$\xi\subset TY$ is a {\sl contact structure} if it is the kernel of a
smooth $1$--form $\th$ on $Y$ such that $\th\wedge d\th \not =0$ at
every point of $Y$\footnote{For an introduction to contact structures
and a guide to the literature we refer the reader to~\cite{Be, El3,
Gi}}. Notice that since $\xi$ is oriented by the restriction of
$d\th$ the manifold $Y$ is necessarily orientable. Moreover, an
orientation on $Y$ induces a coorientation on $\xi$ and
vice-versa. When $Y$ has a prescribed orientation, $\xi$ is said to
be {\em positive} ({\em negative}, respectively), if the orientation on $Y$
induced by $\xi$ coincides with (is the opposite of, respectively) the
given one. In this paper we shall only consider oriented
$3$--manifolds. Therefore, from now on by the expression
``$3$--manifold'' we shall always mean ``oriented $3$--manifold'', and
all contact structures will be implicitly assumed to be positive.

By the work of Martinet and Lutz~\cite{Ma} we know that every closed,
oriented $3$--manifold $Y$ admits a positive contact
structure. Eliashberg defined a special class of contact structures,
which he called {\sl overtwisted}, and proved that in any homotopy
class of cooriented $2$--plane fields on a $3$--manifold there exists
a unique positive overtwisted contact structure up to
isotopy~\cite{El1}. Eliashberg called {\sl tight} the non-overtwisted
contact structures. For tight contact structures, the questions of
existence and uniqueness in a given homotopy class have a negative
answer, in general. For instance, Bennequin proved that there exist
homotopic, non-isomorphic contact structures on $S^3$~\cite{Be},
while Eliashberg showed that the set of Euler classes of tight contact
structures (considered as oriented $2$--plane bundles) on a given
$3$--manifold is finite~\cite{El3}.

The only tight contact structures known at present are fillable in one
sense or another, ie, loosely speaking, they are a $3$--dimensional
phenomenon induced by a $4$--dimensional one. There exist several
different notions of fillability for a contact structure, but here we
shall only define two of them (the weakest ones). The reader
interested in a comprehensive account can look at the
survey~\cite{Et}.

A {\sl $4$--manifold with contact boundary} is a pair $(X,\xi)$, where
$X$ is a connected, oriented smooth $4$--manifold with boundary and
$\xi$ is a contact structure on $\del X$ (positive with respect to the
boundary orientation). A {\sl compatible symplectic form} on $(X,\xi)$
is a symplectic form $\om$ on $X$ such that $\om|_\xi > 0$ at every
point of $\del X$. A contact $3$--manifold $(Y,\ze)$ is called {\sl
symplectically fillable} if there exists a $4$--manifold with contact
boundary $(X,\xi)$ carrying a compatible symplectic form $\om$ and an
orientation-preserving diffeomorphism $\phi$ from $Y$ to $\del X$
such that $\phi_*(\ze)=\xi$. The triple $(X,\xi,\om)$ is said to be a
{\sl symplectic filling} of $Y$. More generally, $(Y,\ze)$ is called
{\sl symplectically semi-fillable} if the diffeomorphism $\phi$ sends
$Y$ onto a connected component of $\del X$. In this case $(X,\xi,\om)$
is called a {\sl symplectic semi-filling} of $Y$. If $(Y,\ze)$ is
symplectically semi-fillable, then $\ze$ is tight by a theorem of
Eliashberg and Gromov (see~\cite{El2,La}).

One of the aims of this paper is to address a fundamental question
about the fillability of contact $3$--manifolds (cf~\cite{El3},
question 8.2.1, and~\cite{Ki}, question 4.142):
\begin{quest}\label{q:fill}
Does every oriented $3$--manifold admit a fillable contact structure? 
\end{quest}
Eliashberg's Legendrian surgery construction~\cite{El1,Go} provides a
rich source of contact $3$--manifolds which are filled by Stein
surfaces (a special kind of $4$--manifolds with contact boundary
carrying exact compatible symplectic forms). Symplectically fillable
contact structures are not necessarily fillable by Stein surfaces. For
example, the $3$--torus $S^1\x S^1\x S^1$ carries infinitely many
isomorphism classes of symplectically fillable contact structures, but
Eliashberg showed~\cite{El4} that only one of them can be filled by a
Stein surface.

Gompf studied
systematically the fillability of Seifert $3$--manifolds using
Eliashberg's construction. This led him to formulate the following:

\begin{conj}[\cite{Go}]\label{conj:Gompf}
The Poincar\'e homology sphere, oriented as the boundary of the
positive $E_8$ plumbing, does not admit positive contact structures
which are fillable by a Stein surface.
\end{conj}

Another basic question asks about the uniqueness of symplectic
fillings. Via Legendrian surgery one can construct, for instance,
non-diffeomorphic (even after blow-up) symplectic fillings of a
given $3$--manifold. On the other hand, $S^3$ is known to have just
one symplectic filling up to blow-ups and
diffeomorphisms~\cite{El2}. We may loosely formulate the uniqueness
question as follows (cf~question 10.2 in~\cite{El2} and question 6
in~\cite{Et}):
\begin{quest}\label{q:uniq}
To what extent does a $3$--manifold determine its symplectic
fillings?
\end{quest}

\subsection{Statement of results}\label{ss:results}

Some progress in the understanding of contact structures has recently
come from studying the spaces of solutions to the Seiberg--Witten
equations. One of the outcomes of~\cite{LM} was a proof of the
existence, for every natural number $n$, of homology $3$--spheres
carrying more than $n$ homotopic, non-isomorphic tight contact
structures. Generalizing to a non-compact setting the results
of~\cite{Ta1, Ta2}, Kronheimer and Mrowka~\cite{KM} introduced
monopole invariants for smooth $4$--manifolds with contact boundary,
and used them to strengthen the results of~\cite{LM} as well as to
prove new results, as for example that on every oriented $3$--manifold
there is only a finite number of homotopy classes of symplectically
semi-fillable contact structures. In this paper we apply~\cite{KM}
to establish the following:

\begin{thm}\label{t:main2}
Let $(X,\xi)$ be a $4$--manifold with contact boundary equipped with a
compatible symplectic form. Suppose that a connected component of the
boundary of $X$ admits a metric with positive scalar curvature. Then,
the boundary of $X$ is connected and $b_2^+(X)=0$.
\end{thm}

The following corollary of theorem~\ref{t:main2} proves
conjecture~\ref{conj:Gompf} as a particular case, and provides a
negative answer to question~\ref{q:fill}.

\begin{cor}\label{c:Poincare}
Let $Y$ denote the Poincar\'e homology sphere oriented as the boundary
of the positive $E_8$ plumbing. Then, $Y$ has no symplectically
semi-fillable contact structures. Moreover, $Y\# -Y$ is not
symplectically semi-fillable with any choice of orientation.
\end{cor}

\begin{proof}
Since $Y$ is the quotient of $S^3$ by a finite group of isometries
acting freely, it has a metric with positive scalar curvature. Hence,
by theorem~\ref{t:main2} if $Y$ is symplectically semi-fillable then
it is symplectically fillable. Moreover, observe that $Y$ cannot be
the oriented boundary of a smooth oriented and negative definite
$4$--manifold. In fact, if $\del X=Y$ then $X\cup (-E_8)$ is a closed,
smooth oriented $4$--manifold with a definite and non-standard
intersection form. The existence of such a $4$--manifold is forbidden
by the well-known theorem of Donaldson~\cite{Do1,Do2}. In view of
theorem~\ref{t:main2}, this proves the first part of the
statement. The second part follows from a general result of
Eliashberg: if $M\# N$ is symplectically semi-fillable, then both $M$
and $N$ are (see~\cite{El2}, theorem~8.1).
\end{proof}

Theorem~\ref{t:main2} can be used, in conjunction with~\cite{Fr}, to
address question~\ref{q:uniq}. Let $(X,\xi)$ be a $4$--manifold with contact
boundary equipped with a compatible symplectic form. Let
$Q_X\co H_2(X;\Z)/\Tor\to\Z$ be the intersection form of $X$. Write the
intersection lattice $J_X=(H_2(X;\Z)/\Tor, Q_X)$ as
\begin{equation*}
J_X=m(-1)\oplus\widetilde{J_X} 
\end{equation*}
for some $m$, where ${\widetilde J_X}$ does not contain classes of
square $-1$. 

\begin{cor}\label{c:fill}
Let $Y$ be a rational homology sphere having a positive scalar
curvature metric. Then, while $X$ ranges over the set of symplectic
fillings of $Y$ such that $\widetilde J_X$ is even, the set of
isomorphism classes of the lattices $\widetilde J_X$ ranges over a
finite set.
\end{cor}

\begin{proof}
By a result of Fr\o yshov (\cite{Fr}, theorem~1) there exists a
rational number $\ga(Y)\in\Q$ depending only on $Y$ such that if $X$
is a negative $4$--manifold bounding $Y$, then for every
characteristic element $\xi\in H_2(X,\del X;\Z)/\Tor$ (ie such that
$\xi\cdot x\equiv x\cdot x\bmod 2$ for every $x\in H_2(X,\Z)/\Tor$),
the following inequality holds:
\begin{equation}\label{e:Froyshov}
\rank(J_X) - |\xi|^2 \leq \ga(Y).
\end{equation}
Thus, if $X$ is a symplectic filling of $Y$, by theorem~\ref{t:main2}
$b_2^+(X)=0$ and therefore equation~\eqref{e:Froyshov}
holds. Clearly~\eqref{e:Froyshov} is also true with $\widetilde J_X$
in place of $J_X$. Hence, if $\widetilde J_X$ is even, choosing
$\xi=0$ we see that the rank of $\widetilde J_X$ is bounded above by a
constant depending only on $Y$. On the other hand, the absolute value
of its determinant is bounded above by the order of $H_1(Y;\Z)$. It
follows (see eg~\cite{MH}) that the isomorphism class of $\widetilde J_X$
must belong to a finite set determined by $Y$.
\end{proof}

\begin{rem}\label{r:oddcase}
{\rm The conclusion of corollary~\ref{c:fill} can be strengthened in
particular cases. For example, if $Y$ is an integral homology sphere,
then the intersection lattice $J_X$ of any symplectic filling of $Y$
is unimodular. It follows from~\cite{Elk1, Elk2} that if $\ga(Y)\leq
8$ then, regardless of whether $\widetilde J_X$ is even or odd, there
are exactly $14$ (explicitly known) possibilities for the isomorphism
class of $\widetilde J_X$ (due to recent work of Mark Gaulter this is
still true as long as $\ga(Y)\leq 24$~\cite{Elk3}). In particular, if
$Y$ is the Poincar\'e $3$--sphere oriented as the boundary of the
negative plumbing $-E_8$, then $\ga(Y)=8$~\cite{Fr}. Up to isomorphism
the only even, negative and unimodular lattices of rank at most eight
are $0$ and $-E_8$. Therefore, $0$ and $-E_8$ are the only
possibilities for $\widetilde J_X$ in this case. Moreover, notice that
if $Y$ bounds a smooth $4$--manifold with $b_2=0$, the same is true
for $-Y$. On the other hand, the argument given to prove
corollary~\ref{c:Poincare} shows that $-Y$ cannot bound negative
semi-definite manifolds. Therefore, if $X$ is an even symplectic
filling of $Y$, $J_X$ is necessarily isomorphic to the negative
lattice $-E_8$.}
\end{rem}

In view of corollary~\ref{c:fill} and remark~\ref{r:oddcase} it seems
natural to formulate the following conjecture:

\begin{conj}\label{conj:mine}
The conclusion of corollary~\ref{c:fill} still holds, under the same
assumptions, if $X$ is allowed to range over the set of all symplectic
fillings of $Y$.
\end{conj}

The plan of the paper is the following. In section~\ref{s:prelim} we
initially fix our notation recalling the results of~\cite{KM}. Then we
state and prove, for later reference, an immediate consequence of
those results, observing how it implies a theorem of Eliashberg. In
section~\ref{s:main} we prove our main result, theorem~\ref{t:main1},
and its corollary theorem~\ref{t:main2}. The line of the argument to
prove theorem~\ref{t:main1} is well-known to the experts. It is the
analogue, in the context of $4$--manifolds with contact boundary, of a
standard argument proving the vanishing of the Seiberg--Witten
invariants of a closed smooth $4$--manifold which splits as a union
$X_1\bigcup_Y X_2$, with $Y$ carrying a positive scalar curvature
metric and $b_2^+(X_i)>0$, $i=1,2$ (cf~\cite{KMT}, remark~6). The
crucial points of such an argument depend on the technical results
of~\cite{MMR}.

{\bf Acknowledgements}. It is a pleasure to thank Dieter Kotschick for
his interest in this paper, and for useful comments on a preliminary
version of it. Warm thanks also go to Peter Kronheimer for observing
that the assumption $b_2^+>0$ in theorem~\ref{t:main1} could be
disposed of when the boundary is disconnected, and to Yasha Eliashberg
for pointing out the second part of
corollary~\ref{c:Poincare}. Finally, I am grateful to the referee for
her/his remarks.

\section{Preliminaries}\label{s:prelim}

We start describing the set-up of~\cite{KM} (the reader is
referred to the original paper for details). A {\sl $\Spin^c$
structure} on a smooth $4$--manifold $X$ is a triple $(W^+,W^-,\ro)$,
where $W^+$ and $W^-$ are hermitian rank--$2$ bundles over $X$ called
respectively the {\sl positive} and {\sl negative spinor bundle}, and
$\ro\co T^* X\to \Hom (W^+,W^-)$ is a linear map satisfying the Clifford
relation: $\ro(\th)^* \ro(\th) = |\th|^2 \Id_{W^+}$ for every $\th\in
T^* X$. The map $\ro$ extends to a linear embedding $\ro\co \La^* T^*
X\to\Hom (W^+,W^-)$. A {\sl $\Spin$ connection} $A$ is a unitary
connection on $W=W^+\oplus W^-$ such that the induced connection on
$\End (W)$ agrees with the Levi--Civita connection on the
image of $\ro$. To any $\Spin$ connection $A$ is associated, via
$\ro$, a twisted Dirac operator $D^+_A\co \Ga(W^+)\to\Ga(W^-)$.

Given a $4$--manifold with contact boundary $(X,\xi)$, let $X^+$ be
the smooth manifold obtained from $X$ by attaching the open cylinder
$[1,+\infty)\x \del X$ along $\del X$. Up to certain choices, the
contact structure $\xi$ determines on $[1, +\infty)\x \del X$ a metric
$g_0$ and a self-dual $2$--form $\om_0$ of constant length
$\sqrt{2}$. $\om_0$ determines on $[1,+\infty)\x \del X$ a $\Spin^c$
structure $\s_0=(W^+,W^-,\ro)$ and a unit section $\Phi_0$ of
$W^+$. Moreover, there is a unique $\Spin$ connection $A_0$ such that
$D^+_{A_0}(\Phi_0)=0$. Given an arbitrary extension of $g_0$ to all of
$X^+$, the triple $(X^+,\om_0,g_0)$ is an AFAK (asymptotically flat
almost K\"ahler) manifold, in the terminology of~\cite{KM}. Consider
the set $\Spin^c(X,\xi)$ of isomorphism classes of $\Spin^c$
structures on $X^+$ whose restriction to $[1,+\infty)\x \del X$ is
isomorphic to $\s_0$. We shall now describe how Kronheimer and Mrowka
define a map
\begin{equation*}
\SW_{(X,\xi)}\co\Spin^c(X,\xi)\to\Z
\end{equation*}
which is an invariant of the pair $(X,\xi)$. Given
$\s=(W^+,W^-,\ro)\in\Spin^c(X,\xi)$, extend $\Phi_0$ and $A_0$
arbitrarily to all of $X^+$. Let $L^2_l$ and $L^2_{l,A_0}$, $l\geq 4$
be, respectively, the standard Sobolev spaces of imaginary $1$--forms
and sections of $W^+$, and let $\mathcal C$ be the space
of pairs $(A,\Phi)$ such that $A-A_0\in L^2_l$ and $\Phi-\Phi_0\in
L^2_{l,A_0}$. Then, ${\mathcal G}=\{u\co X^+\to \C\,|\, |u|=1, 1-u\in
L^2_{l+1}\}$ is a Hilbert Lie group acting freely on ${\mathcal
C}$. Let $\eta\in L^2_{l-1}(i\mathfrak{su}(W^+))$. Given a $\Spin$
connection $A$, let $\hat A$ be the induced $U(1)$ connection on
$\det(W^+)$. Let $M_\eta(\s)$ be the quotient, under the action of
${\mathcal G}$, of the set of pairs $(A,\Phi)\in{\mathcal C}$ which
satisfy the $\eta$--perturbed Seiberg--Witten (or monopole) equations
\begin{equation}\label{e:SW}
\begin{cases}
\ro(F^+_{\hat A})-\{\Phi\otimes\Phi^*\}=
\ro(F^+_{{\hat A}_0})-\{\Phi_0\otimes\Phi_0^*\}+ \eta\\
D_A^+(\Phi)=0,
\end{cases}
\end{equation}
where $\{\Phi\otimes\Phi^*\}$ denotes the traceless part of the
endomorphism $\Phi\otimes\Phi^*$. Kronheimer and Mrowka~\cite{KM}
prove that, for $\eta$ in a Baire set of perturbing terms
exponentially decaying along the end, $M_\eta(\s)$ is (if non-empty) a
smooth, compact orientable manifold of dimension $d(\s)$ equal to
$\langle e(W^+,\Phi_0),[X,\del X]\rangle$, the obstruction to
extending $\Phi_0$ as a nowhere-vanishing section of $W^+$. Now suppose
that an orientation for $M_\eta(\s)$ has been chosen. Then, when
$d(\s)=0$ one can define an integer as the number of points of
$M_\eta(\s)$ counted with signs. $SW_{(X,\xi)}(\s)$ is defined to be
this integer when $d(\s)=0$, and zero when $d(\s)\not=0$.

If $(X,\xi)$ is equipped with a compatible symplectic
form $\om$, then a theorem from \cite{KM} says that there are natural
choices of an element $\s_\om\in\Spin^c(X,\xi)$ and of an orientation
of $M_\eta(\s_\om)$ so that $\SW_{(X,\xi)}(\s_\om)=1$. 

The following proposition is implicitly contained in~\cite{Fr}
and~\cite{KM}. Here we give an explicit statement and proof for the
sake of clarity and later reference.

\begin{prop}\label{p:restr}
Let $(X,\xi)$ be a $4$--manifold with contact boundary. Suppose that
$SW_{(X,\xi)}(\s)\not=0$ for some $\s\in\Spin^c(X,\xi)$. If a
connected component $Y$ of the boundary of $X$ has a metric with
positive scalar curvature then the map $H^2(X;\R)\to H^2(Y;\R)$
induced by the inclusion $Y\subset X$ is the zero map.
\end{prop}

\begin{proof}
The contact structure $\xi$ induces a $\Spin^c$ structure $\t$ on $Y$
(see~\cite{KM}). Let $W$ be the associated spinor bundle on
$Y$. Given a closed $2$--form $\mu$ on $Y$, denote by $N_\mu(Y,\t)$
the set of gauge equivalence classes of solutions to the
$3$--dimensional monopole equations on $Y$ corresponding to the
$\Spin^c$ structure $\t$ and perturbation $\mu$. As observed
in~\cite{KM}, proposition~5.3, it follows from the Weitzenb\"ock
formulae and \cite{Fr} that if $\mu_0\in\Om^2(Y)$ is a closed
$2$--form with $[\mu_0]\not=2\pi c_1(W)$, then there exists a Baire
set of exact $C^r$ forms $\mu_1$ such that $N_{\mu_0+\mu_1}(Y,\t)$
consists of finitely many non-degenerate, irreducible
solutions. Arguing by contradiction, suppose that the restriction map
$H^2(X;\R)\to H^2(Y;\R)$ is non-zero. Then, for every real number
$\ep>0$ there exists a closed $2$--form $\mu$ on $Y$ such that:
\begin{enumerate}
\item[(1)]
$N_\mu(Y,\t)$ consists of finitely many non-degenerate, irreducible
solutions.
\item[(2)]
the $L^2$ norm of $\mu$ is less than $\ep$,
\item[(3)]
$[\mu]\not=2\pi c_1(W)\in H^2(Y;\R)$ and $[\mu]$ is in the image of
the restriction map $H^2(X;\R)\to H^2(Y;\R)$.
\end{enumerate}
Since $\SW_{(X,\xi)}(\s)\not=0$, by~\cite{KM}, proposition~5.8,
$N_\mu(Y,\t)$ is non-empty. But since $Y$ has a metric of
positive scalar curvature, if $\ep$ is sufficiently small the
Weitzenb\"ock formulae imply that $N_\mu(Y,\t)$ is empty: a
contradiction.
\end{proof}

It is interesting to observe that proposition~\ref{p:restr} has the
following corollary, which was first proved by
Eliashberg using the technique of filling by holomorphic disks
\cite{El1}.

\begin{cor}\label{c:S^2xD^2}
$S^2\x D^2$ has no tame almost complex structure with $J$--convex
boundary.
\end{cor}

\begin{proof}
A standard product metric on $S^2\x S^1$ has positive scalar
curvature. Moreover, an almost complex structure on $S^2\x D^2$ has
$J$--convex boundary if, by definition, the distribution $\xi$ of
complex tangents to $S^2\x S^1$ is a positive contact structure. If
$J$ is tame, then there is a compatible symplectic form $\om$ on the
$4$--manifold with contact boundary $(S^2\x D^2,\xi)$. Hence
$\SW_{(S^2\x D^2,\xi)}(\s_\om)\not=0$. But the restriction map
$H^2(S^2\x D^2;\R)\to H^2(S^2\x S^1;\R)$ is non-zero, contradicting
proposition~\ref{p:restr}.
\end{proof}

\section{Proofs of the main results}\label{s:main}

In this section we prove the main results of the paper, namely
theorem~\ref{t:main1} and its immediate corollary,
theorem~\ref{t:main2}. Let $(X,\xi)$ be a $4$--manifold with contact
boundary. We shall start with a preliminary discussion under the
assumption that the boundary of $X$ is connected and admits a metric
with positive scalar curvature. During the proof of
theorem~\ref{t:main1} we will say how to modify the arguments when the
boundary of $X$ is possibly disconnected and at least one of its
connected components admits a metric with positive scalar curvature.

We begin along the lines of~\cite{KM}, proposition~5.6. Let
$(X^+,g_0)$ be the Riemannian $4$--manifold defined in
section~\ref{s:prelim}. We are going to analyze what happens to the
solutions of the equations~\eqref{e:SW} when the metric $g_0$ is
stretched in the direction normal to the boundary of $X$.

In the following discussion we shall denote the boundary of $X$ by
$Y$. Let $g_Y$ be a positive scalar curvature metric on $Y$. Let $g_1$
be a Riemannian metric on $X^+$ coinciding with $g_0$ on
$[1,+\infty)\x Y$ and such that $(X^+,g_1)$ contains an isometric copy
of the cylinder $[-1,1]\x Y$ with the product metric
$dt^2+g_Y$. Choose a perturbing term $\eta_1$ for the monopole
equations which vanishes on this cylinder. For every $R\geq 1$ let
$g_R$ and $\eta_R$ be obtained by replacing $[-1,1]\x Y$ with a
cylinder isometric to $[-R,R]\x Y$. Denote by $X_{\rm in}$ and $X_{\rm
out}$, respectively, the compact and non-compact component of the
complement of the cylinder in $X^+$. Suppose that, for some
$\s\in\Spin^c(X,\xi)$, $\SW_{(X,\xi)}(\s)\not=0$. This implies that
the moduli space $M_{\eta_R}(\s)$ is non-empty for all $R$. Since the
restriction of $\eta_R$ to the cylinder $[-R,R]\x Y$ vanishes, the
proof of lemma~5.7 from~\cite{KM} applies. This says that for every
solution $[A_R,\Phi_R]\in M_{\eta_R}(\s)$ the variation of the
Chern--Simons--Dirac (CSD for short) functional on the restriction of
$[A_R,\Phi_R]$ to $[-R,R]\x Y$ is bounded, independent of $R$. Denote
by $\widehat X_{\rm in}$ and $\widehat X_{\rm out}$ the Riemannian
manifolds obtained by isometrically attaching cylinders $[0,\infty)\x
Y$ and $(-\infty,0]\x \overline Y$ with metric $dt^2+g_Y$ to $X_{\rm
in}$ and $X_{\rm out}$ respectively, where $\overline Y$ denotes $Y$
with the opposite orientation. Let $\eta_{\rm in}$ and $\eta_{\rm
out}$ on $\widehat X_{\rm in}$ and $\widehat X_{\rm out}$ respectively
be compactly supported perturbing terms. Let $R_i$ be a sequence going
to infinity, and let $\eta_i=\eta_{R_i}$ be a corresponding sequence of
perturbing terms as above converging to $\eta_{\rm in}$ and $\eta_{\rm
out}$. Since the moduli spaces $M_{\eta_i}(\s)$ are non-empty for all
$i$, up to passing to a subsequence we may assume that there are
solutions converging on compact subsets to configurations $(A_{\rm
in},\phi_{\rm in})$ and $(A_{\rm out},\phi_{\rm out})$ on $\widehat
X_{\rm in}$ and $\widehat X_{\rm out}$. The configurations $(A_{\rm
in},\phi_{\rm in})$ and $(A_{\rm out},\phi_{\rm out})$ satisfy the
monopole equations for $\Spin^c$ structures $\s_{\rm in}$ and $\s_{\rm
out}$, say, with perturbing terms $\eta_{\rm in}$ and $\eta_{\rm
out}$, and have finite variation of the CSD functional on the
cylindrical ends. Denote the moduli spaces of solutions with bounded
variation of the CSD functional along the end by, respectively,
$M_{\eta_{\rm in}}(\widehat X_{\rm in})$ and $M_{\eta_{\rm
out}}(\widehat X_{\rm out},\xi)$.

The results of~\cite{MMR} imply that $(A_{\rm in},\phi_{\rm in})$,
restricted to the slices $\{t\}\x Y$ converges, as $t\to +\infty$,
towards an element of the moduli space $N_X(Y)$ of solutions of the
unperturbed $3$--dimensional monopole equations on $Y$ modulo the
gauge transformations which extend over $X$. In other words, there is
a map $\del_X\co M_{\eta_{\rm in}}(\widehat X_{\rm in})\to
N_X(Y)$. For every $\th\in N_X(Y)$, we denote $\del_X^{-1}(\th)$ by
$M_{\eta_{\rm in}}(\widehat X_{\rm in},\th)$.

Now recall that, since $\SW(X,\xi)(\s)\not=0$, by the definition of
the invariants $d(\s)=0$, and the canonical spinor $\Phi_0$ can be
extended over $X$ to a nowhere-vanishing section of the bundle
$W^+$. This is equivalent to saying that $\s$ is the $\Spin^c$
structure associated to an almost complex structure $J_X$ on $X$ (see
\cite{KM}, lemma~2.1). Let $Z$ be a smooth, oriented Riemannian
$4$--manifold with boundary $\overline Y$ and such that $J_X$ extends
to an almost complex structure $J_M$ on the closed oriented
$4$--manifold $M=X\cup_Y Z$ (the reason why such a $Z$ exists is
explained in eg~\cite{Go}, lemma 4.4; one can always find a $Z$ such
that the obstruction to extending $J_X$ over $Z$ is concentrated at a
finite number of points, and then, in order to kill the obstruction,
one can modify $Z$ by connect summing at those points with a suitable
number of copies of $S^2\times S^2$). Let $\widehat Z$ be the manifold
with cylindrical end obtained by attaching $(-\infty, 1]\x\overline Y$
to the boundary of $Z$. Fix an extension of $J_M$ from $Z$ to
$\widehat Z$, and call $\s_{\widehat Z}$ the $\Spin^c$ structure
induced on $\widehat Z$. Choose an identification of the cylindrical
ends of $\widehat X_{\rm out}$ and $\widehat Z$ (observe that
$\s_{\widehat Z}$ is isomorphic to $\s_{\rm out}$ on the cylindrical
end). Also, choose a perturbing term $\eta'$ on $\widehat Z$ which
coincides with $\eta_{\rm out}$ on the cylindrical end. As before,
there is a moduli space $M_{\eta'}(\widehat Z)$, a map $\del_X\co
M_{\eta'}(\widehat Z)\to N_Z(\overline Y)$, and, for every $\th'\in
N_Z(\overline Y)$, we denote $\del_Z^{-1}(\th')$ by
$M_{\eta'}({\widehat Z},\th')$.

\begin{lem}\label{l:dimsum}
For any $\th_1\in N_X(Y)$, $\th_2\in N_Z(\overline Y)$, $M_{\eta_{\rm
in}}(\widehat X_{\rm in},\th_1)$ and $M_{\eta'}({\widehat Z},\th_2)$
are (possibly empty) smooth manifolds. Moreover, the sum of their
expected dimensions equals $-1-b_1(Y)$.
\end{lem}

\begin{proof}
By a standard argument (see eg~\cite{MMS}), since the metric $g_Y$
has nowhere negative scalar curvature, the moduli space $N_X(Y)$
consists of reducible solutions, and the linearization of the
equations on $Y$ with appropriate gauge fixing gives a deformation
complex whose first cohomology group at a point $[A,0]\in N_X(Y)$ can
be identified with $H^1(Y;\R)\oplus\ker D_A$. Since $g_Y$ has positive
scalar curvature, we have $\ker D_A=0$ for every $[A,0]\in
N_X(Y)$. Moreover, since the dimension of $N_X(Y)$ is $b_1(Y)$,
$N_X(Y)$ is smooth, and the Kuranishi map from the first to the second
cohomology of the deformation complex vanishes. It follows
from~\cite{MMR} that every element of $M_{\eta_{\rm in}}(\widehat
X_{\rm in})$ converges, along the end, exponentially fast towards an
element of $N_X(Y)$. This implies that, given any $\th\in N_X(Y)$,
$M_{\eta_{\rm in}}(\widehat X_{\rm in},\th)$ is a (possibly empty)
smooth manifold. Exactly the same arguments apply to
$M_{\eta'}({\widehat Z})$.

Recall that taking the quotient of $N_X(Y)$ by the whole gauge group
of $Y$ gives a covering map $p\co N_X(Y)\to N(Y)$ with fiber
$H^1(Y;\Z)/H^1(X;\Z)$. For every $\th_1\in N_X(Y)$, denote $p(\th_1)$
by ${\overline\th_1}$. Let $W_X^+$ be the spinor bundle associated
with the $\Spin^c$ structure $\s_{\rm in}$. By~\cite{APS}
and~\cite{MMR} the exponential convergence implies that, given
$\th_1=[A,0]$, the expected dimension of $M_{\eta_{\rm in}}(\widehat
X_{\rm in},\th_1)$ is
\begin{equation}\label{e:dim1}
d_1=\frac14 (c_1(W_X^+)^2-2\chi(X)-3\si(X))-\frac{h^0({\overline\th_1})+ 
h^1({\overline\th_1})}2 + \eta_Y({\overline\th_1})
\end{equation}
where $h^0({\overline\th_1})=1$ is the dimension of the stabilizer of
the configuration $(A,0)$, and $h^1({\overline\th_1})=b_1(Y)$ is the
dimension of the first cohomology group of the deformation complex at
$(A,0)$. $\eta_Y({\overline\th_1})$ is the $\eta$--invariant of the
relevant boundary operator on $Y$ defining the deformation complex
(since we are going to use only well known properties of this
operator, we don't need to be more specific, see~\cite{MMS} for more
details). Note that the rational number $c_1(W_X^+)^2$ is well defined
because by proposition~\ref{p:restr} $c_1(W_X^+)|_Y$ is a torsion
class.

Similarly, if $\th_2\in N_Z({\overline Y})$, the expected dimension of
$M_{\eta'}({\widehat Z},\th_2)$ is
\begin{equation}\label{e:dim2}
d_2=\frac14 (c_1(W_Z^+)^2-2\chi(Z)-3\si(Z))-\frac{h^0({\overline\th_2})+ 
h^1({\overline\th_2})}2 + \eta_{\overline Y}({\overline\th_2}).
\end{equation}
Again, $h^0({\overline\th_2})=1$ and
$h^1({\overline\th_2})=b_1(Y)$. Recall that $\eta_Y$ changes sign when
the orientation of $Y$ is reversed. Moreover, since
$h^0({\overline\th})$ and $h^1({\overline\th})$ are constant in
${\overline\th}\in N(Y)$ there is no spectral flow, and therefore
$\eta_Y(\overline\th)$ is constant too. Hence, $\eta_{\overline
Y}({\overline\th_2})=-\eta_Y({\overline\th_2})=-\eta_Y({\overline\th_1})$.
Finally, observe that the $\Spin^c$ structures $\s_{\rm in}$ and
$\s_Z$ can be glued together to give a $\Spin^c$ structure $\s_M$ on
the closed manifold $M=X\cup_Y Z$. In fact, $\s_M$ can be taken to be
the $\Spin^c$ structure induced by the almost complex structure $J_M$
(see the discussion before the statement). It follows
that the associated spinor bundle $W_M^+$ satisfies
\begin{equation*}
c_1(W_M^+)^2=2\chi(M)+3\si(M),
\end{equation*}
and the formula $d_1+d_2=-1-b_1(Y)$ follows immediately
from~\eqref{e:dim1} and~\eqref{e:dim2}.
\end{proof}

\begin{thm}\label{t:main1}
Let $(X,\xi)$ be a $4$--manifold with contact boundary. Suppose that 
one of the following assumptions holds:
\begin{enumerate}
\item[\bf1\rm)] The boundary of $X$ is connected, it admits a
metric with positive scalar curvature and $b_2^+(X)>0$,
\item[\bf2\rm)]
The boundary of $X$ is disconnected and one of its connected
components admits a metric with positive scalar curvature.
\end{enumerate}
Then, the map $\SW_{(X,\xi)}$ is identically zero.
\end{thm}

\begin{proof}
We will start by establishing the conclusion under the first
assumption. Arguing by contradiction, suppose that the map
$\SW_{(X,\xi)}$ does not vanish. Then, one can argue as in~\cite{KM},
proposition~5.4, and show that, for $\eta_{\rm in}$ in a Baire set of
compactly supported perturbations, if, for some $\th_1\in N_X(Y)$,
$M_{\eta_{\rm in}}(\widehat X_{\rm in},\th_1)$ is non-empty, then its
expected dimension is non-negative (observe that, since the
perturbing term is decaying to zero along the cylindrical end, we need
$b_2^+ (X)>0$ to rule out reducible solutions). Thus, choosing
$\eta_{\rm in}$ in such a Baire set, the existence of $(A_{\rm
in},\Phi_{\rm in})$ implies $d_1\geq 0$. If we denote by $d_2$ the
expected dimension of $M_{\eta_{\rm out}}(\widehat X_{\rm
out},\xi,\th_2)$ (with the obvious meaning of the symbols), the same
argument gives $d_2\geq 0$ (no assumption on $b_2^+$ is needed now,
because the elements of $M_{\eta_{\rm out}}(\widehat X_{\rm
out},\xi,\th_2)$ are asymptotically irreducible on the ``conical''
end). As explained in~\cite{KM}, subsection~5.4, one can associate to
$\th_2$ a homotopy class of $2$--plane fields $I(\th_2)$ on $Y$. As in
the proof of proposition~5.6 in~\cite{KM}, the expected dimension of
$M_{\eta_{\rm out}}(\widehat X_{\rm out},\xi,\th_2)$ is given by a
difference element ${\overline\de}\left(I(\th_2),\xi\right)$
(see~\cite{KM}, subsection~5.1, for the definition of
$\overline\delta$; in the case at hand this number is an integer
because, by proposition~\ref{p:restr}, the restriction of $c_1(W^+)$
to $Y$ is a torsion element). Moreover,
${\overline\de}\left(I(\th_2),\xi\right)$ is also equal to the
expected dimension of $M_{\eta'}({\widehat Z},\th_2)$. This
contradicts lemma~\ref{l:dimsum}. Hence, we have established the
conclusion of the theorem under the first assumption.

When the boundary of $X$ is disconnected the above argument can be
easily modified so that the requirement on $b_2^+(X)$ becomes
redundant. In fact, one can repeat the same construction involving
only the end corresponding to the boundary component having positive
scalar curvature. $\widehat X_{\rm in}$ will have one cylindrical end
as well as some conical ends $E_i$, $i=1,\ldots,k$, while $\widehat
X_{\rm out}$ will be the same as before. The conical ends can be
chopped off and replaced by suitable compact manifolds with boundary
$Z_i$ (as we did before with $\widehat X_{\rm out}$) without changing
the expected dimension of the corresponding moduli spaces. Then,
denoting $\left(\widehat X_{\rm in}\setminus\cup E_i\right)\cup Z_i$
by $\widetilde X_{\rm in}$, the statement of lemma~\ref{l:dimsum} will
still hold with $M_{\eta_{\rm in}}\left(\widehat X_{\rm
in},\th_1\right)$ replaced by $M_{\eta_{\rm in}}\left(\widetilde
X_{\rm in}, \th_1\right)$, and will have a similar proof. On the other
hand, the same arguments as before show that, for generic choices of
$\eta_{\rm in}$, the expected dimensions of $M_{\eta_{\rm
in}}\left(\widehat X_{\rm in},\xi_1,\ldots,\xi_k,\th_1\right)$ (with
the obvious meaning of the symbols) and $M_{\eta_{\rm out}}(\widehat
X_{\rm out},\th_2,\xi)$ are non-negative, and they coincide with the
expected dimensions of $M_{\eta_{\rm in}}\left(\widetilde X_{\rm in},
\th_1\right)$ and $M_{\eta'}({\widehat Z},\th_2)$, respectively. No
assumption on $b_2^+(X)$ is needed, because both $\widehat X_{\rm in}$
and $\widehat X_{\rm out}$ have at least one conical end, and the
elements of $M_{\eta_{\rm in}}\left(\widehat X_{\rm
in},\xi_1,\ldots,\xi_k,\th_1\right)$ and $M_{\eta_{\rm out}}(\widehat
X_{\rm out},\xi,\th_2)$ are asymptotically irreducible on the conical
ends. This gives a contradiction as in the previous case, and
concludes the proof of the theorem.
\end{proof}

\begin{proof}[Proof of theorem~\ref{t:main2}]
Let $\om$ be the compatible symplectic form. We know (see
section~\ref{s:prelim}) that there is a distinguished element
$\s_\om\in\Spin^c(X,\xi)$ such that $\SW_{(X,\xi)}(\s_\om)\not=0$. The
conclusion follows immediately from theorem~\ref{t:main1}.
\end{proof}

%
%


\begin{thebibliography}

\bibitem{APS} {\bf M\,F Atiyah}, {\bf V\,K Patodi}, {\bf I\,M Singer}, {\it
Spectral asymmetry and Riemannian geometry: I}, Math. Proc. Cambridge
Philos. Soc. 77 (1975) 43--69

\bibitem{Be} {\bf D Bennequin}, {\it Entrelacements et equations de Pfaff},
Ast\'erisque 107--108 (1983), 83--161

\bibitem{Do1} {\bf S\,K Donaldson}, {\it Connections, cohomology and
the intersection forms of four--manifolds}, Jour. Diff. Geom. 24
(1986) 275--341

\bibitem{Do2} {\bf S\,K Donaldson}, {\it The Seiberg--Witten equations
and $4$--manifold topology}, Bull. AMS 33 (1996) 45--70

\bibitem{El1} {\bf Y Eliashberg}, {\it Topological characterization of 
Stein manifolds of dimension $>2$}, Intern. Journal of Math. 1, 
No. 1 (1990) 29--46

\bibitem{El2} {\bf Y Eliashberg}, {\it Filling by holomorphic discs and its
applications}, London Math. Soc. Lecture Notes Series 151 (1991) 45--67

\bibitem{El3} {\bf Y Eliashberg}, {\it Contact $3$--manifolds twenty
years since J.~Martinet's work}, Ann. Inst. Fourier 42 (1992) 165--192

\bibitem{El4} {\bf Y Eliashberg}, {\it Unique holomorphically fillable
contact structure on the $3$--torus}, Intern. Math. Res. Not. 2 (1996)
77--82

\bibitem{Elk1} {\bf N Elkies}, {\it A characterization of the $\Z^n$
lattice}, Math. Res. Lett. 2 (1995) 321--326

\bibitem{Elk2} {\bf N Elkies}, {\it Lattices and codes with long
shadows}, Math. Res. Lett. 2 (1995) 643--652

\bibitem{Elk3} {\bf N Elkies}, personal communication

\bibitem{Et} {\bf J\,B Etnyre}, {\it Symplectic convexity in low
dimensional topology}, Top. Appl. (to appear)

\bibitem{Fr} {\bf K\,A Fr\o yshov}, {\it The Seiberg--Witten equations and
four--manifolds with boundary}, Math. Res. Lett. 3 (1996)
no. 3, 373--390

\bibitem{Gi} {\bf E Giroux}, {\it Topologie de contact en dimension
$3$}, S\'eminaire Bourbaki 760 (1992-93), 7--33

\bibitem{Go} {\bf R\,E Gompf}, {\it Handlebody construction of Stein
surfaces}, Ann. of Math. (to appear)

\bibitem{Ki} {\bf R Kirby}, {\it Problems in Low-Dimensional
Topology}. In W H Kazez (Ed.), Geometric Topology, Proc of the 1993
Georgia International Topology Conference, AMS/IP Studies in Advanced
Mathematics, pp.~35--473, AMS \& International Press
(1997)

\bibitem{KM} {\bf P\,B Kronheimer}, {\bf T\,S Mrowka}, {\it Monopoles and
contact structures}, Invent.~Math. 130 (1997) 209--256

\bibitem{KMT} {\bf D Kotschick}, {\bf J\,W Morgan}, {\bf C\,H Taubes}, {\it
Four--manifolds without symplectic structures but with non-trivial
Seiberg--Witten invariants}, Math. Res. Lett. 2 (1995) 119--124

\bibitem{La} {\bf F Laudenbach}, {\it Orbites p\'eriodiques et courbes
pseudo-holomorphes, application \`a la conjecture de Weinstein en
dimension $3$ [d'apr\`es H.~Hofer et al.]}, Ast\'erisque 227
(1995) 309--333

\bibitem{LM} {\bf P Lisca}, {\bf G Mati\'c}, {\it Tight contact structures
and Seiberg--Witten invariants}, Invent. math. 129 (1997) 509--525

\bibitem{Ma} {\bf J Martinet}, {\it Formes de contact sur les 
vari\`et\`es de dimension $3$}, Lect. Notes in Math. 209, 
Springer--Verlag (1971) 142--163

\bibitem{MH} {\bf J Milnor}, {\bf D Husemoller}, {\it Symmetric bilinear
forms}, Ergebnisse der Mathematik und Ihrer Grenzgebiete, Band 73,
Springer--Verlag (1973)

\bibitem{MMR} {\bf J\,W Morgan}, {\bf T\,S Mrowka}, {\bf D
Ruberman}, {\it The $L^2$--moduli space and a vanishing theorem for
Donaldson polynomial invariants}, Monographs in Geometry and Topology,
no. II, International Press, Cambridge, MA, 1994

\bibitem{MMS} {\bf J\,W Morgan}, {\bf T\,S Mrowka}, {\bf Z
Szab\'o}, {\it Product formulas along $T^3$ for Seiberg--Witten
invariants}, preprint (1997)

\bibitem{Ta1} {\bf C\,H Taubes}, {\it The Seiberg--Witten invariants and 
symplectic forms}, Math. Res. Lett. 1 (1995) 809--822

\bibitem{Ta2} {\bf C\,H Taubes}, {\it More constraints on symplectic
manifolds from Seiberg--Witten equations}, Math. Res. Lett. 2
(1995) 9--14

\end{thebibliography}
\end{document}